\documentstyle[11pt]{article}
\newtheorem{Co}{Corollary}
\newtheorem{Prop}{Proposition}
\newtheorem{Df}{Definition}
\newtheorem{Th}{Theorem}
\newtheorem{Le}{Lemma}

\def\ov{\overline}
\newtheorem{Ex}{Example}
\def\bc{\begin{centre}}
\def\ec{\end{centre}}
\def\bq{\begin{equation}}
\def\eq{\end{equation}}
\def\cd{\cdot}
\def\ba#1{\begin{array}{#1}}
\def\ea{\end{array}}

\def\b{\beta}
\def\a{\alpha}
\def\g{\gamma}

\def\s{\sigma}

\def\p{\varphi}

\def\ps{\psi}

\def\l{\lambda}
\def\r{\rho}

\def\c{\cal}

\def\t{\cal T}
\def\li{\cal L}
\def\m{\cal M}

\def\G{\cal G}
\def\ov{\overline}
\def\sm{\setminus}
\def\t{\bigtriangleup}
\def\op{\oplus}
\def\es{\emptyset}
\def\suo#1{\sum_{#1\in\t}\op}
\def\bo{\Box}
\def\su#1{\sum_{#1\in\t}}
\def\ot{\otimes}
\def\L{\Lambda}
\def\ra{\rightarrow}

\begin{document}

%\newpage
\begin{center}
{\LARGE {\bf Lie algebras with triality.}}
\end{center}

\vspace{2cm}

\centerline{\large {\bf Alexander Grishkov}\footnote{supported fy FAPESP and
CNPq (Brazil)}}

\begin{center}
Rua do Mat\~ao,1010,IME,
 University of S\~ao Paulo, CEP 005598-300,
S\~ao Paulo, Brazil.
\end{center}
\begin{center}
e-mail: grishkov@ime.usp.br
\end{center}
\begin{center}
{\it In memorian Valerii Filippov}
\end{center}

\newpage

\section{ABSTRACT}

By analogy with the definition of group with triality we introduce
Lie algebra with triality as Lie algebra  $L$ wich admits the group of automorphisms
$S_3=\{\s,\,\r |\s^2=\r^3=1,\,\s\r\s=\r^2\}$ such that for any
$x\in L$ we have $(x^{\s}-x)+(x^{\s}-x)^{\rho}+(x^{\s}-x)^{\rho^2}=0.$
We describe the structure of finite dimensional Lie algebra with triality
over a field of characteristic 0 and give applications of Lie algebras with triality to the
theory of Malcev algebras.
\newpage
\section{Introduction.}
Glauberman \cite{Gl1}, \cite{Gl2} noted that the multiplication group $G(M)$ of a given Moufang loop
$M$ admits a certain dihedral group $S_3$ of automorphisms,
with $|S_3|=6.$ For fixed generators $\s,\r\in S_3$ with
$\s^2=\r^3=1,\, \s\r\s=\r^2,$ the following equation holds for any $x$ in $G(M):$
\bq
(x^{\s}x^{-1})(x^{\s}x^{-1})^{\r}(x^{\s}x^{-1})^{\r^2}=1.
\label{per}
\eq
S.Doro \cite{Dor} called groups with such automorphisms {\em groups with triality},
since the most striking example is $D_4(q)$ with its graph automorphisms.

In this paper we introduce the notion of a Lie algebra with triality
 such that the simple Lie algebra of type $D_4$ is a Lie algebra with triality and
a Lie algebra that corresponds to an algebraic or Lie group
with triality is a Lie algebra with triality. By definition,
a Lie algebra $L$ over a field of characteristic $\neq 2,3$ is a Lie algebra with triality
if $L$ admits the dihedral group $S_3$ of automorphisms such that the following
analog of (\ref{per}) holds for any $x$ in $L:$

\bq
(x^{\s}-x)+(x^{\s}-x)^{\rho}+(x^{\s}-x)^{\rho^2}=0.
\label{per2}
\eq
The Lie algebras with this property first appeared in \cite{Mih}
where a connection of those algebras with Malcev algebras was described.
For any Lie algebra $L$ we can construct the following Lie algebra with triality
$T(L)=L_1\oplus L_2\oplus L_3,$ where $L_i$ is an isomorphic copy
of $L$ with fixed isomorphism $l\ra l_i,\,l\in L,\,i=1,2,3.$
The algebra $T(L)$ admits the following action of $S_3:$
$$
l_1^{\s}=l_2,\,l_2^{\s}=l_1,\,l_3^{\s}=l_3,\,l_1^{\rho}=l_2,\,l_2^{\rho}=l_3,\,l_3^{\rho}=l_1.
$$
It is clear that $T(L)$ with this $S_3$-action is a Lie algebra with triality. We call a
Lie algebra with triality $A$
$standard$ if it is isomorphic to either some algebra with triality of type $T(L)$ or its invariant
subalgebra.
One of the main results of this paper is the following Theorem
\begin{Th}
Let $L$ be a perfect finite dimensional
Lie algebra with triality  over an algebraicly closed field of characteristic 0.
Then $L$ is an extension of a Lie algebra with triality of type
$D_4$ and a standard Lie algebra with triality.
\end{Th}
For more exact formulation see Theorem 5.

In the third section we prove that for a Lie algebra with triality $L$
the condition (\ref{per2}) is equivalent to the condition
$L=L_0\oplus L_2$ where $L_0=\{x\in L | x^{\l}=x,\,\forall \l\in S_3\}$ and
$L_2$ is a sum of irreducible 2-dimensional $S_3$-modules. We call this decomposition
an $S_3$-$decomposition.$
Hence the notion of Lie algebra with triality is a particular case of the following
general definition.
\begin{Df}
Let $S$ be a finite group, $\{\L_1, ... ,\L_m\}$
be a set of non-isomorphic absolutly irreducible $S$-modules, and assume $\L=\L_1\oplus ... \oplus \L_m.$
An algebra $A$ is called a $\L$-algebra if $S\subseteq AutA$ and the $S$-module $A$
has a decomposition $A=A_1\oplus ... \oplus A_n$ such that each
$S$-module $A_i$ is isomorphic to one of the $S$-modules in $\{\L_1, ... ,\L_m\}.$
\end{Df}
In the second section of this paper we develop a method of studying
the $\L$-algebras belonging to a given variety ${\cal N}$.
The main idea is the following: we fix a finite group $S,$ a set of irreducible
$S$-modules $\{\L_1, ... ,\L_m\},$ and a variety ${\cal N}.$
For any $\L$-algebra $A$ one can construct an algebra $H_{\L}(A)=Hom_S(\L,A)$
(see Proposition 2 for a detailed definition) which has the following grading
$H_{\L}(A)=\sum_{i=1}^m\oplus Hom_S(\L_i,A).$
We denote by ${\cal N}_{\L}$ the category of all $\L$-algebras
from ${\cal N}$ and by $H_{\L}({\cal N})$ the category of corresponding graded
algebras $H_{\L}({\cal N})=\{H_{\L}(A) | A\in {\cal N}_{\L}\}.$
These categories ${\cal N}_{\L}$ and $H_{\L}({\cal N})$ are isomorphic
but the advantage of $H_{\L}({\cal N})$ is that it admits a structure of
graded variety.

In the third section we apply the method of studying the $\L$-algebras
defined in the second section to the theory of Lie algebras with triality.
The main result of the third section is the construction of a functor
$\Psi$ from the category ${\cal L}$ of Lie algebra with triality to the category ${\cal T}$
of graded algebras and a functor $\Psi_0$ from the category ${\cal T}$
to the category ${\cal M}$ of Malcev algebras
(see (\ref{7}) for definition of Malcev algebra).
These functors allow us to obtain the main results of the theory
of finite dimensional Malcev algebras over a field of characteristic 0
as corollaries of the
corresponding well known results on
finite dimensional Lie algebras (see Theorems 1 and 2).

In the fourth section we apply the correspondence between
Lie algebra with triality and Malcev algebras to the theory of
Lie algebra with triality.
We introduce a notion of $T$-centre of Lie algebra with triality
which is an analogue of the notion of Lie centre in the theory of
Malcev algebras. Observe that the $T$-centre of a given Lie algebra with triality
is an invariant ideal and it is a trivial algebra with triality.
In Theorem 4 we prove that every perfect finite dimensional
Lie algebra with triality over a field of characteristic 0 without
$T$-centre is a Lie algebra with triality of type $D_4.$
In the last part of this section we apply deep results
of the theory of infinite dimensional Malcev algebras obtained by  V. Filippov \cite{Fi1}-\cite{Fi8}
to the theory of Lie algebra with triality. In particular, we prove that
every simple (non necessarily finite dimensional) Lie algebra with triality over a field
of characteristic $\neq 2,3$ is a Lie algebra with triality of type $D_4.$

All spaces and algebras are considered over a fixed field $k$ of
characteristic
$\neq 2,3.$ Further,
$k\{X\}=kX$ denotes the $k-$space with
a basis $X.$

\section{Graduate variety}

In this section we prove some results about
 graduate varieties and give their applications
 to
Lie algebras.

We fix a set $\t$ and call any
space $V$ with a fixed $\t$-grading:
$V=\sum_{\a\in \t}\op V_{\a}$ $\t$-space (resp. $\t$-algebra, $\t$-module).
We can consider the $\t$-space $V$ as an algebra with unary operations
$\{\a |\a\in\t\}$ such that for $a\in V$ $(a)_{\a}=a_{\a}$ if
$a=\su\b a_{\b}.$
 Let $A=\sum_{\a\in\t}\op A_{\a}$
be a $\t$-algebra. Then a $\t$-identity for the $\t$-algebra $A$ is
a (non-associative) polynomial $f(x, y, ...)$
in signature $(+,\cd,\a\in \t)$
such that
$f(a, b, ...)=0,$ for all elements
$a,b, ... \in A.$ For example, $f(x,y)=(x_{\a}y_{\b})_{\g}+(x_{\tau}^2)_{\g}.$
Let $V=\suo{\a} V_{\a}, W=\suo\b W_{\b}$ be given $\t$-spaces.
We define the $contraction$ of the $\t$-spaces as follows:
$$V\bo W=\sum_{\a\in \t}\oplus (V_{\a}\otimes W_{\a}).$$
Thus  $V\bo W$ is a $\t$-space too.

If $A$ and $B$ are two $\t$-algebras then we define the
$contraction$ of $\t$-algebras $A$ and $B$ as a $\t$-space $A\bo B$ with the following
multiplication rule:
$$(a_{\a}\ot b_{\a})\cd (a_{\b}\ot b_{\b})=
\sum_{\g\in\t} c_{\g}\ot d_{\g},$$
where
$a_{\a}a_{\b}=\su\g c_{\g},\, b_{\a}b_{\b}=\su\g d_{\g}.$
\begin{Df}
A set $\c N$ of algebras over $k$ is called a $\t$-variety if
$\c N$ is the set of all $\t$-algebras over $k$ which satisfy
a given set of $\t$-identities.

For a given set $X$ of $\t$-algebras or $\t$-identities
we denote by $\{X\}$ the minimal
$\t$-variety that contains $X$ or
satisfies all identities from $X$.

If $\c N$ and $\c M$ are two $\t$-varieties then we define
a contraction operation by
$${\c N} \bo {\c M }=\{A\bo B | A\in{\c N} , B\in{\c M}\},$$
and a division operation by
\[
\ba{ll}
{\c N} /{\c M} =  & \{ A  |\ \  \forall B\in{\c M} , A\bo B\in{\c N}\ \
              \mbox{and}\ \  A \ \ \mbox{satisfies} \\
& \mbox{all identities of} \ \
{\c M} \ \ \mbox{of the type} \ \ (a_{\a}b_{\b})_{\g}=0\}.
\ea
\]
\end{Df}
It is obvious that
$(\c N /\c M)\bo\c M\subset \c N.$

 Between these two operations (contraction and division) there is some
difference. If we have a set $X$ of $\t$-polynomials
such that $\c N $=$\{X\}$ and a $\t$-algebra $A$ such that
$\c M$ =$\{A\}$ then finding a set $Z$ such that $\c N\bo\c M$ =$\{Z\}$ may be non-trivial.
 On the other hand there is a simple algorithm
for constructing a set of $\t$-identities $Z$ such that
$\c N /\c M$ =$\{Z\}$.
 First we have to take the absolutely free
$\t$-algebra $F=F(x_1, ...),$ where $\{x_i\}$ are the  homogeneous free generators
of $F$.
Let $B=\{a_i| i=1,...\}$ be a homogeneous basis of the $\t$-algebra $A$.
For any $\t$-identity $f(x_{\a_1}, ... , x_{\a_n})$ of the set $X$
and any subset $T=\{a_{i_1}, ... , a_{i_n}\}$ of $B$ such that
$a_i\in A_{\a_i}$ we construct the following set of $\t$-identities:
$G(f,T)=\{g_1(x_{\a}, ...), ... ,g_m(x_{\a}, ... )\}$,
where
$f(x_{\a_1}\ot a_{i_1}, ... ,x_{\a_n}\ot a_{i_n})=
\sum_{j=1}^{m}g_j\ot a_{k_j}$ and $k_j\neq k_i,$ if $j\neq i$.

\begin{Prop}
If $\c N$ and $\c M$ are $\t$-varieties such that
$\c N $=$\{X\}$, $\c M =$$\{A\}$ and $B=\{a_i|i=1, ... \}$
 is a homogeneous basis of the $\t$-algebra $A$  then
$$
{\c N}/ {\c M}=\{M_2, G(f,T)\,|\,f\in X, T\subset B\},
$$
where $M_2$
is the set of identities of the variety $\c M$ of the type
$(x_{\a}x_{\b})_{\g}=0$.
\end {Prop}
{\bf Proof.}
If $C\in {\c N}/{\c M}$ then $D=C\bo A\in {\c N}$ and, by definition,
$C$ satisfies all identities from $G(f,T)$ for $f\in X.$
Conversely, if we have a $\t$-algebra $C$ which satisfies all identities
from $G(f,T)$ for $f\in X$ then the $\t$-algebra $C\bo A$ satisfies
the identity $f.$ Hence $C\bo A\in {\c N}.$
$\Box$

We note that if $\c N$ is some variety (not necessarily
graded) then, by definition, a $Z_2$-graded algebra $A$ is a
$\c N$-superalgebra if the algebra $A\bo G\in \c N$, where
$G$ is a Grassmannian algebra. That is the $Z_2$-variety
${\c N}_2$ of $\c N$-superalgebras is ${\c N}/{\c G}$, where
${\c G}=\{G\}$. It is well known that
 there is an easy algorithm to construct the graded identities of
the $\c N$-superalgebras if we know the identities of the variety
$\c N$.

Fix a finite group $S$ such that $char(k)=p$ is not a divisor of $|S|$ and some set
of absolutly irreducible non-isomorphic $S$-modules $\L(S)=\{\L_1, ... ,\L_m\}.$ We denote
$\L=\L_1\oplus ... \oplus \L_m$ and say that an $S$-module
$V$ is of type $\L$ if $V=V_1\oplus ... \oplus V_n,$
$V_1, ... ,V_n\in \L(S)$
and of type $\ov{\L}$ if $V_i\notin \L(S), i=1, ... n.$
It is clear that for any $S$-module $V$ we have
$V=V_{\L}\op V_{\ov{\L}},$ where $V_{\L}$ is the submodule in $V$
of type $\L$ and $V_{\ov{\L}}$ is one of type $\ov{\L}.$
 We write $V\otimes_{\L}W$ for $(V\otimes W)_{\L}.$

Suppose that the set $\L(S)$ has the following property:
for every $i,j,k\in \{1, ... ,m\}$ $dimHom_S(\L_i\otimes_{\L}\L_{j},\L_k)<2$
then
there exists an embedding
\bq
\phi_{ij}: \L_i\otimes_{\L}\L_{j}\ra\L.
\label{1}
\eq
It is clear that $\L=\L_{ij}\oplus \ov{\L_{ij}}$ where
$\L_{ij}=\phi_{ij}(\L_i\otimes_{\L}\L_{j}).$
We define
$$\psi_{ij}:\L\ra \L_i\otimes_{\L}\L_{j}$$
by $\psi_{ij}(\ov{\L_{ij}})=0$ and
$\psi_{ij}(\phi_{ij}(v))=v,\,v\in \L_i\otimes_{\L}\L_{j}.$
Here and above we identify $v\otimes w$ (for $v\in V,w\in W$)
with its image in $V\otimes_{\L}W.$
The module $\L$ has the following structure of algebra

\bq
v\cd w=\phi(v\otimes w)=\phi_{ij}(v\otimes w), \,\,v\in \L_i, w\in \L_j.
\label{2}
\eq

On the other hand, $\L$ has the following co-algebraic structure:

\bq
\phi^{\star}(v)=\sum_{ij}\psi_{ij}(v).
\label{3}
\eq

\begin{Df}
An algebra $A$ is called an algebra of type $\L$ if
$S\subseteq Aut_kA$ and $A$ is an $S$-module of type $\L.$
\end{Df}
We denote by ${\c N}_{\L}$ the category of all algebras
of type $\L$ from a given variety ${\c N}.$
If $\t=\{1, ... ,m\}$ then all algebras of ${\c N}_{\L}$
are $\t$-algebras and morphisms of the category ${\c N}_{\L}$
preserve this graduation.
In general the category ${\c N}_{\L}$ is not a $\t$-variety.
This means that set of $\t$-identities
$X$ such that ${\c N}_{\L}=\{X\}$ there is no.
But there exists a natural isomorphism (as categories!) between ${\c N}_{\L}$
and some $\t$-variety.
\begin{Prop}
Let ${\c N}$ be a variety and ${\c N}_{\t}$ be the $\t$-variety of
all $\t$-algebras from ${\c N}.$
Then the map
$\Psi: {\c N}_{\t}/\{\L\}\ra {\c N}_{\L}$
defined by $\Psi(C)=C\bo \L$ is an isomorphism between
the categories ${\c N}_{\t}/\{\L\}$ and ${\c N}_{\L}.$
\end{Prop}
{\bf Proof.}
It is clear that $\Psi$ is a functor from ${\c N}_{\t}/\{\L\}$
into ${\c N}_{\L}.$
To finish the proof of the Proposition, it is enough
to construct the inverse functor
$\Phi:{\c N}_{\L}\ra {\c N}_{\t}/\{\L\}.$ For any
$A\in {\c N}_{\L},$ we define $\Phi(A)=Hom_S(\L,A)$
and for $\xi,\tau\in Hom_{S}(\L,A)$ we define
$\xi\star\tau\in Hom_S(\L,A)$
so, that the following diagram is commutative
\bq
\ba{ccc}
\L&
\stackrel{\phi^{\star}}{\longrightarrow}&\L\otimes_{\L}\L\\
\Big\downarrow
\vcenter{%
\rlap{$\xi\star\tau$}}& &
\Big\downarrow
\vcenter{\rlap{$\xi\otimes\tau$}}\\
A &
\stackrel{m}{\longleftarrow} & A\otimes_{\L}A\\
\ea
\label{4}
\eq
Here $m:A\otimes_{\L}A\ra A $ is the multiplication in $A$ and
$\phi^{\star}$ is the comultiplication in $\L.$

We need to prove that for $A\in {\c N}_{\L}$ we have
$\Psi(\Phi(A))\simeq A.$ Define a linear map $\pi:\Psi(\Phi(A))\ra A$
as follows
$\pi(\xi_i\otimes v_i)=\xi_i(v_i)$ if
$\xi_i\otimes v_i\in \Psi(\Phi(A))=\Phi(A)\bo\L=Hom_S(\L,A)\bo\L,$
$\xi_i\in Hom_S(\L_i,A),$ $v_i\in \L_i.$
As $A$ is an algebra of type $\L$ it is clear that $Im(\pi)=A.$
Let $A=\sum_{j=1}^m\op A_j$ where $A_j$ is a sum of $S$-modules of type $\L_j.$
Hence
\[
\ba{l}
dim_kHom_S(\L_j,A)=dim_kHom_S(\L_j,A_j)=\\[4mm]
dim_kA_j/dim_k\L_j=k_j.
\ea
\]
Therefor
\[
\ba{l}
dim_k\Psi(\Phi(A))=dim_kHom_S(\L,A)\bo\L=
\sum_{j=1}^mdim_k(Hom_S(\L_j,A)\otimes\L_j)=\\[2mm]

\sum_{j=1}^mk_jdim_k\L_j=
\sum_{j=1}^mdim_kA_j=dim_kA.
\ea
\]
Hence $ker(\pi)=0$ and $\pi$ is a linear isomorphism of linear spaces.
We next prove that $\pi$ is an isomorphism of algebras.
Choose $\xi_i\in Hom_S(\L_i,A),\xi_j\in Hom_S(\L_j,A),$
$v_i\in \L_i,v_j\in\L_j$ and suppose that
$\xi_i\star\xi_j=\sum_{s=1}^m\l_s,$
$v_i\cd v_j=\sum_{s=1}^mw_s,$
where $\l_s\in Hom_S(\L_s,A), w_s\in \L_s,s=1, ... ,m.$

By definition, we have
$$
(\xi_i\otimes v_i)\cd(\xi_j\otimes v_j)=\sum_{s=1}^m\l_s\otimes w_s,
$$
hence
\bq
\pi(\xi_i\otimes v_i)\cd\pi(\xi_j\otimes v_j)=\xi_i(v_i)\xi_j(v_j).
\label{5}
\eq
On the other hand,
$\pi((\xi_i\otimes v_i)\cd(\xi_j\otimes v_j))=\pi(\sum_{s=1}^m\l_s\otimes w_s)=\sum_{s=1}^m\l_s(w_s).$
 From (\ref{4}) we have
$(\xi_i\star\xi_j)(v_i\cd v_j)=\xi_i(u_i)\cd\xi_j(u_j)$
if $\phi^{\star}(v_i\cd v_j)=\sum_{p,q}u_p\otimes u_q,$
but  from (\ref{3}) we get

$$\phi^{\star}(v_i\cd v_j)=\phi^{\star}(\phi(v_i\otimes v_j))=v_i\otimes v_j+r_{ij},$$
where $r_{ij}\in \sum_{(p,q)\neq (i,j)}\L_p\otimes_{\L}\L_q.$
Hence
$$\pi((\xi_i\otimes v_i)\cd(\xi_j\otimes v_j))=(\xi_i\star\xi_j)(v_i\cd v_j)=\xi_i(v_i)\xi_j(v_j).$$
From this and (\ref{5}) we have that $\pi$ is a homomorphism.
$\Box$

{\bf Note.} It is clear that all constructions and Propositions above are true
if we substitute a finite group $S$ by an algebra $S$ such that all finite dimensional
$S$-modules are semisimple. For example, instead of a finite group $S$ we can take
any finite dimensional
semisimple Lie, Malcev or Jordan algebra.

\begin{Ex}
Let $S$ be a 3-dimensional simple Lie algebra with a basis $\{e,h,f\}$ and
multiplication law $eh=2e,ef=h,fh=-2f.$
We fix an $S$-module
\linebreak
$\L=S\op V,$ as above, where
$V$ has a basis $\{v_1,v_{-1}\}$ with $S$-action $v_1e=v_{-1}f=0,
v_1f=v_{-1},v_{-1}e=-v_1,v_1h=v_1,v_{-1}h=v_{-1}.$
Then $\L$ has a structure of algebra of type $\L$
such that $S$ is a subalgebra, $V$ is an $S$-submodule with the action
as above and
$$
vx=-xv,\forall x\in S,v\in V,
$$
$$
v_1^2=e,v_1v_{-1}=v_{-1}v_1=h,v_{-1}^2=f.
$$
Note that $\L$ has a ${\bf Z}_2$-graduation $\L=\L_0\op\L_1,$
where $\L_0=S$ and $\L_1=V$ and this algebra is a Lie superalgebra
$osp(1,2).$
\end{Ex}
The following Proposition is an easy corollary of the
Propositions 1 and 2.
\begin{Prop}
A Lie algebra $L$ is an algebra of type $\L=S\op V$ if and only if
$L=H(L)\bo\L,$ where $H(L)\in {\c L}/\{\L\}={\cal G}r$ and ${\cal G}r$
is a ${\bf Z}_2$-variety of Grassmannian algebras.
\end{Prop}

\begin{Ex}
Let $S$ be the 3-dimensional simple Lie algebra as in Example 1
and $\Gamma=S\op W$ be an $S$-module such that
$W$ is the unique irreducible Malcev $S$-module
with a basis $\{w_2,w_{-2}\}$ and
 with $S$-action $w_{-2}e=w_{2}f=0,
w_2e=-2w_{2},w_{-2}f=2w_2,w_2h=2w_2,w_{-2}h=-2w_{-2}.$
Then $\Gamma$ has a structure of algebra of type $\Gamma$
such that $S$ is a subalgebra, $W$ is an $S$-submodule with the action
as above and
$$
vx=-xv,\forall x\in S,v\in W,
$$
$$
w_2^2=f,w_2w_{-2}=w_{-2}w_2=h,w_{-2}^2=e.
$$
Note that $\Gamma$ has a ${\bf Z}_2$-graduation $\Gamma=\Gamma_0\op\Gamma_1,$
where $\Gamma_0=S$ and $\Gamma_1=W.$
\end{Ex}
The following analogue of Proposition 3 was proved in \cite{Gr5}.
\begin{Prop}
A Malcev algebra $M$ is an algebra of type $\Gamma=S\op W$ if and only if
$M=H(M)\bo\Gamma,$ where $H(M)\in {\c M}/\{\L\}={\c NG}$ and ${\c NG}$
is a ${\bf Z}_2$-variety
with the following ${\bf Z}_2$-identities
\bq
\ba{l}
an=na,xy=-yx,\\[2mm]
(an)m=a(nm),\\[2mm]
(xy)z+(yz)x+(zx)y=0,
\ea
\label{6}
\eq
where $a\in H(M)_0,x,y,z\in H(M)_1,n,m\in H(M).$
\end{Prop}

\section{ Lie algebra with triality.}

In this section we will consider the central example of application of Propositions 1 and 2.
  Recall that an algebra $A$ is an algebra with triality
if this algebra admits a nontrivial action of the group
$S_3=\{\s,\r \,|\,\s^2=\r^3=1, \s\r\s=\r^{-1}\}$
by automorphisms such that for every $x\in A$ we have (\ref{per2}).
In what follows, fix a set $\t=\{0,2\}$ and genarators $\s,\rho$ of $S_3.$

\begin{Le}
Let $A$ be an algebra over a field of characteristic $p\neq 2,3$
and $S_3\subseteq Aut_kA.$
Then $A$ is a algebra with triality if and only if
$A=A_0\oplus A_2,$ where $A_0=A^{S_3}=\{a\in A\,|\, a^g=a,\forall g\in S_3\}$
and $A_2=\sum\oplus V_i,$ $V_i$ is an irreducible  two dimensional $S_3$-module.
\end{Le}
{\bf Proof.}
Let $V$ be $S_3$-module with standard basis $\{v,w\}$ such that
$$
v^{\s}=w, w^{\s}=v, v^{\r}=w, w^{\r}=-v-w.
$$
Then $V$ is the unique irreducible 2-dimensional $S_3$-module.
Let $A=A_0\oplus A_2$ and $x\in A.$ Then there exist irreducible $S_3$-modules
$V_1$ and $V_2$
with standard bases
$x_v,x_w$ and $y_v,y_w,$
respectively, such that $x=x_0+x_v+y_w,$
where $x_0\in A_0.$
Then $x^{\s}-x=x_w-x_v-(y_w-y_v)=z_w-z_v,$ where $z_w=x_w-y_w,\,z_v=x_v-y_v.$
 Hence

 $(x^{\s}-x)+(x^{\s}-x)^{\rho}+(x^{\s}-x)^{\rho^2}=z_w-z_v-z_v-z_w-z_w+2z_v+2z_w=0.$

Conversely, let $A$ be a algebra with triality and let $x\in A$ be such that
$x^{\s}=-x,\, x^{\r}=x.$ Then from (\ref{per2}) we have

$0=(x^{\s}-x)+(x^{\s}-x)^{\rho}+(x^{\s}-x)^{\rho^2}=-6x.$

Hence $x=0.$ $\Box$

 Fix a $\t$-graded algebra $\L=\L_0\oplus \L_2,$ where $\L_2=V$ is the irreducible $S_3$-module
with standard basis $\{v,w\}$
as above, $\L_0=ka$
and
$$
a^2=a, av=va=v,aw=wa=w,
$$
\bq
v^2=(v+2a)/3, w^2=(w+2a)/3, vw=wv=-(v+w+a)/3.
\label{f1}
\eq
$\L$ is a algebra with triality with this $S_3$-action.

It is clear that an algebra $A$ is a algebra with triality if and only if it is an algebra
of type $\L=\L_0\op \L_2.$
For  any Lie algebra with triality $A=A_0\op A_2$ we denote
$N(A)=Ann_{A_0}A_2$ and $K(A)=A_2\oplus (A_2^2)_0/N(A)\cap (A_2^2)_0.$ It is clear that $N(A)$ is an ideal of $A$
and $N(A/N(A))=0.$
We call a Lie algebra with triality $A$ $normal$ if $N(A)=0$ and $A_0=(A_2^2)_0.$
It is obvious that $A$ is normal if and only if $A=K(A).$

We give some examples of Lie algebra with triality.

\begin{Ex}
Let $L$ be a Lie algebra and $L_1,L_2,L_3$ be Lie algebras
isomorphic to $L.$
As in the introduction $A=T(L)=L_1\oplus L_2\oplus L_3$ has the structure
of Lie algebra with triality.
In this case
\[
\ba{l}
A_0=\{l_1+l_2+l_3 | l\in L\},\\[3 mm]

A_2=\{a_1+b_2+c_3 |a,b,c\in L, a+b+c=0\},\\[3 mm]

N(A)=\{l_1+l_2+l_3 | l\in Z(L)=Ann_LL\}.
\ea
\]
Observe that $A$ is normal if and only if $Z(L)=0$ and $L^2=L.$
\end{Ex}
\begin{Df}
A Lie algebra with triality $P$ is called {\bf trivial} if there exists a Lie algebra $L$
such that $K(P)\simeq K(B)$ where $B$ is an invariant subalgebra of $T(L).$
\end{Df}
 To describe the set of invariant subalgebras
of $T(L)$ we need the following definition.
\begin{Df}
Let $L$ be a Lie algebra. Then a pair of subalgebras $(A,B)$ of $L$ is {\bf compatible}
if $B^2\subseteq A$ and $AB\subseteq B.$
\end{Df}
 For any compatible pair $(A,B)$ of subalgebras in a Lie algebra $L$ we
construct an invariant subalgebra $T(A,B)\subseteq T(L).$
By definition, $T(A,B)=T(A,B)_0\oplus T(A,B)_2$ where
\[
\ba{l}
T(A,B)_0=\{l_1+l_2+l_3 \,|\,l\in A \},
\ea
\]
and
\[
\ba{l}
T(A,B)_2=\{a_1+b_2+c_3 |a,b,c\in B,\,a+b+c=0 \}.
\ea
\]
\begin{Prop}
Let $L$ be a Lie algebra and $T(L)$ be the corresponding standard Lie algebra with triality.
Then for every invariant subalgebra $P$ in $T(L)$
there exists a unique compatible pair $(A,B)$ of subalgebras in $L$
such that $P=T(A,B)$
\end{Prop}
{\bf Proof.}
Let $P=P_0\oplus P_2$ be the $S_3$-decomposition of $P.$
We denote $A=\{l\in L\,|\,l_1+l_2+l_3\in P_0\}$
and $B=\{l\in L\,|\,l_1-l_2\in P_2\}.$
It is obvious that $A$ is a subalgebra of $L$ and $AB\subseteq B.$
Let $l,r\in B.$ We have
\[
\ba{l}
(l_1-l_2)(r_1-r_2)=(lr)_1+(lr)_2=\\[3mm]
2[(lr)_1+(lr)_2+(lr)_3]/3+[(le)_1+(lr)_2-2(lr)_3]/3\in P.
\ea
\]
Hence $lr\in A\cap B$ and $(A,B)$ is a compatible pair.$\Box$

As a corollary of this Proposition we construct examples of
trivial Lie algebra with triality $P$ and $Q$ such that
$P$ and $Q/N(Q)$ are standard Lie algebra with triality but
$P/N(P)$ and $Q$ are not standard. Recall that a Lie algebra with triality is standard
if it is an invariant subalgebra of a Lie algebra with triality of type $T(L).$
\begin{Ex}
Let $L=sl_p(k)$ be a Lie algebra of matrices
with trace zero over a field of characteristic $p>3.$
Then the centre $Z(L)$ of $L$ is the unique proper ideal
of $L$ and $Z(L)=ke.$ We define $P=T(L)$ and $Q=P/I$
where $I=k(e_1-e_2)\oplus k(e_1-e_3).$
If $Q$ is a standard Lie algebra with triality then by Proposition 5
there exists a compatible pair $(A,B)$ of some Lie algebra $R$
such that
$Q=T(A,B).$ If $C=A\cap B\neq 0$ then $Q$ has an invariant
ideal $J=T(C)$ and $dim_kZ(J)=3dim_kZ(C).$
On the other hand, $Q$ has only two non-zero invariant ideals
$k(e_1+e_2+e_3)$ and $Q,$ which have one dimensional centre.
This is a contradiction.

Hence $A\cap B=0$ and $AB=0.$ Then
$Q_2=\{\a l_1+\b l_2+\g l_3 \,|\,l\in B,\a,\b,\g\in k, \a+\b+\g=0 \}$ is an invariant
ideal, which is a contradiction. Thus we proved that $Q$ is not a standard
Lie algebra with triality. It is obvious that $Q/N(Q)=T(L/Z(L))$ is a standard
Lie algebra with triality.

Analogously we can prove that $P/N(P)$ is not a standard Lie algebra with triality.
\end{Ex}

\begin{Ex}
Let $L$ be a split simple finite dimensional Lie algebra of type $D_4.$
In \cite{G8} the following basis of $L$ over a field of characteristic $\neq 2$ was constructed

$$B=\{e_i,h_i,f_i,i=1,...,4;\,\mu |\mu\subseteq I_4=(1234)\}$$
with the following multiplication law for the basis elements
\[
\begin{array}{l}
e_i\cd f_i=h_i,\, e_i\cd h_i=2e_i,\,h_i\cd f_i=2f_i,\\[3mm]
e_i\cd \p= \p\cup i,\,\,i\in I_4\sm\p;\\[3mm]
\p\cd f_i=\p\sm i,\,\,i\in \p;\\[3mm]
\p\cd h_i=\p,\,\,i\in \p;\\[3mm]
\p\cd h_i=-\p\,\,,i\in I_4\sm \p;
\end{array}
\]
\[
\p\cd \ps=\left\{
\begin{array}{ll}
(-1)^{|\ps|+1}e_i,& \p\cap\psi=i,\p\cup\ps=I_4;\\[3mm]
(-1)^{|\ps|}f_i,& \p\cap\psi=\es,\p\cup\ps=I_4\sm i;\\[3mm]
(-1)^{|\ps|}(\sum_{i\in\ps}h_i-\sum_{j\in \p}h_j)/2,& \p\cap\psi=\es,\p\cup\ps=I_4.
\end{array}
\right.
\]
Here $|\s|$ is the number of elements of $\s\subseteq \{1,2,3,4\}.$

The natural action of $S_3$ on $I_4=\{1,2,3,4\}$
such that $1^{\s}=2,\,1^{\r}=2,\,2^{\r}=3,\,3^{\r}=1,$
$4^{\l}=4,\,\forall \l\in S_3$ can be extended to the set
$\{\mu | \mu\subseteq I_4\}.$
We can extend this $S_3$-action to $B$ so, that $x_i^{\l}=x_{\l(i)},$
$\l\in S_3,\,i\in I_4,\,x\in \{e,f,h\}.$
It is easy to see that each element of $S_3$ acts on $L$ as an automorphism.
We introduce an order on $B_0=B\sm \{h_1, ... ,h_4\}$ such that
$b>0$ if and only if $b\in \{e_1, ... ,e_4\}$ or $b\in \{\mu | \mu\subseteq I_4,\,4\in\mu\}.$
If we identify $B_0$ with the set of roots of $L$ then $B_+=\{b\in B_0 | b>0\}$ is the set
of the positive roots and $\{4,e_1,e_2,e_3\}$ is the set of an simple roots.
Since $e_i\cd e_j=0$ and $e_i\cd 4=\{i4\},i\neq 4$ the corresponding Dynkin
diagram is of type $D_4$ and $S_3$ acts as diagram automorphisms.
Applying Lemma 1 it is easy to prove that $L$ is an algebra with triality.
\end{Ex}

By Proposition 2 we have.
\begin{Prop}
An algebra $A$ (not neassarily Lie) is an algebra with triality if and only if
there exists a $\t$-graded algebra $M=M_0\oplus M_2,$ such that
$M_0^2\subset M_0,\,\,M_0M_2\subset M_2,\,\,M_2M_0\subset M_2,$\,\,
$A=M\bo \L$ and $S_3$ acts on the first term of the product.
\end{Prop}

 By Proposition 1, the set of all $\t$-graded algebras $M=M_0\oplus M_2$
such that $M_0^2\subset M_0,\,\,M_0M_2\subset M_2, \,\,M_2M_0\subset M_2$ and
$M\bo \L$ is a Lie algebra, forms a $\t$-graded variety ${\c T}$.
\begin{Prop}
A $\t$-graded algebra $M=M_0\op M_2\in {\c T}$ if and only if
$M$ satisfies the following $\t$-graded identities:
$$(ab)_2=(xa)_0=(ax)_0=0,$$
\bq
m^2=0,
\label{f2}
\eq
\bq
(mn)_ia=(ma\cd n)_i+(m\cd na)_i,\,\,i=0,2,
\label{f3}
\eq
\bq
6(xy)_0z=((xy)_2z)_2+((zy)_2x)_2+((xz)_2y)_2,
\label{f4}
\eq
\bq
((xy)_2z)_0+((zx)_2y)_0+((yz)_2x)_0=0,
\label{f5}
\eq
where $n,m\in M, a\in M_0, x,y,z\in M_2$ and $n=(n)_0+(n)_2,$ for
$(n)_0\in M_0, (n)_2\in M_2.$
\end{Prop}
{\bf Proof.}
Let $C\in {\c T};$ then $C=C_0\op C_2$ and $C\bo\L=C_0\otimes \L_0\op C_2\otimes \L_2$
is a Lie algebra. Then for any $x,y,z\in C_2$ and $b,c\in C_0$ we have
$$
[[c\ot a,x\ot v],y\ot w]+[[x\ot v,y\ot w],c\ot a]+[[y\ot w,c\ot a],x\ot v]=0.
$$
From this and (\ref{f1}) we have
\[
\ba{l}
[cx\ot v,y\ot w]-[(xy)_0\ot a,c\ot a]/3-[(xy)_2\ot(v+w),c\ot a]+\\[2mm]

[yc\ot w,x\ot v]=-(((cx)y)_0\ot a)/3-(((cx)y)_2\ot(v+w))/3-\\[2mm]

((xy)_0c\ot a)/3-((xy)_2c\ot(v+w))/3-(((yc)x)_0\ot a)/3-\\[2mm]

(((yc)x)_2\ot(v+w))/3=0.
\ea
\]
Hence
$$
(xy)_0c+((cx)y)_0+((yc)x)_0=0,
$$
and
$$
(xy)_2c+((cx)y)_2+((yc)x)_2=0.
$$
Analogously, from
$$
[[z\ot w,x\ot v],y\ot v]+[[x\ot v,y\ot v],z\ot w]+[[y\ot v,z\ot w],x\ot v]=0
$$
and (\ref{f1}), we have

\[
\ba{l}
-\{((xy)_2z)_0+((yz)_2x)_0+((zx)_2y)_0\ot a\}/9+\\[4mm]
\{2(xy)_0z/3-((xy)_2z)_2/9+((yz)_2x)_2/9+((zx)_2y)_2/9\}\ot w+\\[4mm]
\{-((xy)_2z)_2/9-(yz)_0x/3-(zx)_0y/3\}\ot v=0.

\ea
\]
Then the $\t$-identities (\ref{f4}) and (\ref{f5}) hold. Note that we have one more identity
\[
((xy)_2z)_2=3(yz)_0x/3-3(zx)_0y.
\]
But this identity is equivalent to (\ref{f4}).

It is easy to see that the same calculations prove that an algebra
$C\bo \L$ is a Lie algebra if $C$ satisfies the $\t$-identities (\ref{f2}-\ref{f5}).
$\Box$

  We denote by ${\c L}={\cal L}(k)$  the category of Lie algebra with triality over a field $k$
and by ${\c M}={\cal M}(k)$ the category of Malcev algebras
over $k.$
From the Propositions 2 and 6 we have functors
$\Psi:\li\ra{\c T}$ and $\Phi=\Psi^{-1}:{\c T}\ra\li,$ where by definition,
$\Phi(A)=A\bo\L.$

The following Lemma gives a new proof of Theorem 1 \cite{Mih}.
\begin{Le}
If $M=M_0\oplus M_2\in {\cal T}$ then the space $M_2$ with the product
$x\star y=(xy)_2$ is a Malcev algebra.
\end{Le}
{\bf Proof.}
Recall that an anti-commutative algebra $A$ is a Malcev algebra
if it satisfies the identity

\bq
((xy)z)x+((yz)x)x+((zx)x)y=(xz)(yx).
\label{7}
\eq

Multiplying (\ref{f5}) by $6x$ and applying (\ref{f4}) we obtain:
\bq
\ba{l}
6((xy)_2z)_0x+6((yz)_2x)_0x+6((zx)_2y)_0x=\\[3mm]
(((xy)_2z)_2x)_2+((xz)_2(xy)_2)_2+(((xy)_2z)_2x)_2+\\[3mm]
2(((yz)_2x)_2x)_2+(((zx)_2y)_2x)_2+(((zx)_2x)_2y)_2+\\[3mm]
((xy)_2(zx)_2)_2=\\[3mm]
((x\star y)\star z)\star x+2((y\star z)\star x)\star x+
((x\star y)\star x)\star z+\\[3mm]
((z\star x)\star y)\star x+2(x\star y)\star (z\star x)+
((z\star x)\star x)\star y=0.
\ea
\label{g6}
\eq
Note that, in general, an anticommutative algebra with identity (\ref{g6})
is not a Malcev algebra. Thus we need some more identities.

By (\ref{f4}) we have
\bq
3(xy)_0x=((xy)_2x)_2,  \ \ 3(xy)_0y=((xy)_2y)_2.
\label{f7}
\eq
Now from (\ref{f3}) and (\ref{f7}) we can obtain

\bq
\ba{lll}
3(xy)_0(xy)_2=&3(((xy)_0x)y)_2+ &3(x\cd (xy)_0y)=\\[3mm]
((((xy)_2x)_2)y)_2+&(x\cd ((xy)_2y)_2)_2.&
\ea
\label{g8}
\eq
On the other hand, (\ref{f4}) yields
$$
6(xy)_0(xy)_2=(((xy)_2y)_2x)_2+((x\cd (xy)_2)_2y)_2.
$$
 By this and (\ref{g8}) we have
\bq
3((x\star y)\star y)\star x=3((x\star y)\star x)\star y.
\label{f9}
\eq
Let $J(x,y,z)=(x\star y)\star z+(y\star z)\star x+(z\star x)\star y.$
Then (\ref{g6}) and (\ref{f9}) can be rewritten as
\bq
2J(x,y,z)\star x+J(x\star y,x,z)+J(z\star x,x,y)=0,
\label{f10}
\eq
\bq
J(x\star y,x,y)=0.
\label{f11}
\eq
By the linearization of (\ref{f11}) we have
$$
J(x\star y,x,z)+J(x\star z,x,y)=0.
$$
Then by (\ref{f10}) one obtains
\bq
J(x,y,z)\star x=J(x\star z,x,y).
\label{f12}
\eq
But the identities (\ref{7}) and (\ref{f12}) are equivalent.
Note that identities (\ref{f10}) and (\ref{f11}) follow from (\ref{f12}).
$\Box$

Using Lemma 2 we define a functor $\Psi_0: {\cal T}\ra {\cal M}.$
If $M=M_0\op M_2\in {\cal T}$ let $\Psi_0(M)=(M_2,\star).$
We denote by ${\c F}$ the functor $\Psi_0\circ\Psi:{\cal L}\ra{\cal M}.$
Note that this functor was constructed in another way by Mikheev \cite{Mih}.

 Now we construct the left inverse functor $\G:\m\ra\li$ of ${\c F}.$
For this we define a functor $\Phi_0:{\cal M}\ra {\cal T}$ so, that, for $M_2\in {\cal M}$,
$\Phi_0(M_2)=Inder(M_2)\op M_2.$
Here $Inder(M_2)$ is the Lie algebra of the inner derivations
of $M_2$:
$$Inder(M_2)=\{D(x,y)=L_{[x,y]}+[L_x,L_y]\| x,y\in M_2,L_x:y\ra xy\}.$$
The multiplication law $\cd$ in $\Phi_0(M_2)$ is defined by the standard action of
$Inder(M_2)$ on $M_2$ and the following formula:
\bq
x\cd y=D(x,y)/6+[x,y].
\label{f6}
\eq
\begin{Prop}
$\Phi_0$ is a functor from $\m$ into ${\c T}.$
\end{Prop}
{\bf Proof.}
By definition, if $M_2\in \m,$ then the algebra
$M=\Phi_0(M_2)=Inder(M_2)\op M_2=M_0\op M_2$
has a $\t$-gradation with the following $\t$-identities:
$(ab)_2=(xa)_0=0$ for $a,b\in M_0,x\in M_2.$
Hence it is enough to prove that $M$ satisfies the $\t$-identities
(\ref{f2})-(\ref{f5}).
The $\t$-identities (\ref{f2}), (\ref{f3}) and (\ref{f4}) hold in $M$ by definition.
The identity (\ref{f5}) can be rewritten in the form:
$$
D([x,y],z)+D([y,z],x)+D([z,x],y)=0,x,y,z\in M_2,
$$
or
\bq
\ba{l}
[x,y,z,t]+[[t,z],[x,y]]+[x,y,t,z]+[z,x,y,t]+[[t,y],[z,x]]+\\[3mm]
[z,x,t,y]+[y,z,x,t]+[[t,x],[y,z]]+[y,z,t,x]=0.
\ea
\label{f8}
\eq
But (\ref{f8}) is a consequence of the complete
linearization
of the identities (\ref{7}) and (\ref{f11}).
$\Box$

Finalty, we define ${\cal G}=\Phi\circ\Phi_0.$
It is clear that ${\cal F}({\cal G}(M_2))=M_2$ if $M_2\in \m$ but, in general,
${\cal G}({\cal F}(L))\not=L,$ for $L\in \li.$

\begin{Prop}
Let $L\in \li.$ Then ${\cal G}({\cal F}(L))=L$ if and only if
$L$ is normal.
\end{Prop}
{\bf Proof.}
Let $L$ be a normal Lie algebra with triality. Then $L=M\bo \L$ where
$M=M_0\op M_2\in {\c T}.$ As $L$ is normal we have that
$Ann_{M_0}M_2=0$ and $M_0=(M_2^2)_0.$
From (\ref{f4}), we have that for $x,y\in M_2$ the operator
$L_{(xy)_0}: M_2\ra M_2$ is equal to the inner derivation
$D(x,y)/6.$ As $M_0=(M_2^2)_0,$ we have a homomorphism
$\pi:M_0\ra Inder(M_2),$
$\pi(\sum (x_iy_i)_0)=\sum D(x_i,y_i).$
It is clear that $ker(\pi)=Ann_{M_0}M_2=0.$
Hence $M=\Phi_0(M_2)$ and ${\c G}({\c F}(L))={\c G}(M_2)=L.$

The converse is obvious.
$\Box$

Note that in general the variety ${\cal T}$ does not lie in ${\cal L}.$
 However we have the following
Proposition.

\begin{Prop}
Let $P$ be a Lie algebra with triality over a field $k$ of characteristic either 0 or a prime $p>5.$
Let $M= \Psi(P)\in {\cal T}$ and $M_2=\Psi_0(M)\in {\cal M}(k).$ Then the following
conditions are equivalent:

(i) The Malcev algebra $M_2$ is a Lie algebra;

(ii) $P$ is a trivial Lie algebra with triality;

(iii) M is a Lie algebra.

\end{Prop}

{\bf Proof.}
Let $M_0\oplus M_2\in {\cal T}$ then $M$ is a Lie algebra if and only if
for every $x,y,z\in M_2$ we have
\bq
\ba{ll}
xy\cd z+yz\cd x+zx\cd y=&
((xy)_2z)_2+((yz)_2x)_2+((zx)_2y)_2+\\[4mm]
(xy)_0z+(yz)_0x+(zx)_0y=&0.
\ea
\label{t1}
\eq
Multiplying (\ref{t1}) by 6 and applying (12)
we can obtain:
\bq
\ba{l}
6((xy)_2z)_2+6((yz)_2x)_2+6(((zx)_2y)_2)+((xy)_2z)_2+\\[2mm]
((zy)_2x)_2+((xz)_2y)_2+((yz)_2x)_2+((xz)_2y)_2+\\[2mm]
((yx)_2z)_2+((zx)_2y)_2+((yx)_2z)_2+((zy)_2x)_2=\\[2mm]
5((xy)_2z)_2+5((yz)_2x)_2+5(((zx)_2y)_2)=0.
\ea
\label{x}
\eq
Since the field $k$ has characteristic 0 or a prime $p>5,$ we have that $M$ is a Lie algebra if
and only if $M_2$ is a Lie algebra.

Now suppose that $M$ and $M_2$ are Lie algebras. We need to prove that $P$ is a
trivial Lie algebra with triality. We can suppose that $K(P)=P.$
It is clear that $((M_2^2)_2,M_2)$ is a compatible pair of $M_2.$
Let $A$ be the corresponding invariant subalgebra of $T(M_2).$
Then $K(A)\simeq {\cal G}(M_2)\simeq K(P).$ Hence $P$ is a trivial
Lie algebra with triality.
$\Box $

As a consequence of equality (\ref{x}) we have
\begin{Co}
Let $k$ be a field of characteristic 5 and $P$ be a Lie algebra with triality over $k.$
Then $\Psi(P)$ is a Lie algebra and ${\cal F}(P)$ is a Lie algebra
if and only if $P$ is a trivial Lie algebra with triality.
\end{Co}

Let $k$ be a field of characteristic 0 or $p>3$ that contains an element
$\xi\neq 1$ such that $\xi^3=1.$
Then for any Lie algebra over $k$ and $\r\in Aut_{k}L,$ $\r^3=1,$ we have
a ${\bf Z}_3-$gradation of $L:$ $L=L_1\oplus L_0\oplus L_{-1}$ where
$L_i=\{a\in L | a^{\r}=\xi^ia\}.$ We will say that $\r$ admits triality if there exists an involution
$\s\in Aut_kL$ such that $\s\r\s=\r^2$ and $a^{\s}=a,\,\forall a\in L_0.$ We note that in this case
$L$ is a Lie algebra with triality with respect to action of $S_3=\{\r,\s\}.$ Note that in this case
$L_1^{\s}=L_{-1},$ $L_{-1}^{\s}=L_{1}.$ Indeed, if $x\in L_1$ then
$x^{\s\r}=x^{\s\r\s\s}=x^{\r^2\s}=\xi^2x^{\s}.$ By Proposition 7 we can write $L=M\bo \L,$ for some
$M\in {\cal T}.$ Then $L_0=M_0\otimes \L_0,$ $L_1=M_2\otimes k(v-\xi w),$
$L_{-1}=M_2\otimes k(v-\xi^2 w).$ It is clear that $L_1$ and $L_{-1}$ as $L_0-$modules are isomorphic.

Suppose that we have an automorphism of order 3 $\r$ of $L$ such that  $L_1$ and $L_{-1}$ as
$L_0-$modules are isomorphic. Let us fix an $L_0-$isomorphism $\tau:L_1\to L_{-1}.$ In this case
the $k-$vector space $L_1$ has two algebra structures: $(L_1,\circ)$ and $(L_1,\bullet),$ where,
by definition, for $x,y\in L_1:$ $x\circ y=[x,y]^{\tau},$ $x\bullet y=[x^{\tau},y^{\tau}].$

\begin{Le}
In notations above $\rho$ admits triality if and only if there exists an $L_0-$isomorphism
$\tau:L_1\to L_{-1}$ such that $x\circ y=x\bullet y,\,[x,x^{\tau}]=0,\,\forall x,y\in L_1.$
\end{Le}
{\bf Proof.}
Suppose that $\r$ admits triality and $\s$ is a corresponding involution.
Then $L=M\bo \L,$
$L_i=M\otimes k(v-\xi^iw),i=1,2,$ and $\s:L_1\to L_{-1},\,\s(m\otimes (v-\xi w))=-\xi m\otimes (v-\xi^2w),$
is an $L_0-$isomorphism. For $x,y\in L_1$ we get $x\circ y=[x,y]^{\s}=[x^{\s},y^{\s}]=x\bullet y. $
It is clear that for $x=m\otimes (v-\xi w)\in L_1$ we get
$[x,x^{\s}]=0$ since $m^2=0.$

Let $\r,\tau$ satisfy the hypothesis of Lemma. We define an involution $\s$ of $L$ by the following:
$x^{\s}=x,\,x\in L_0,$ $x^{\s}=x^{\tau},\,(x^{\tau})^{\s}=x,x\in L_1.$ We have,
since $\tau$ is $L_0-$isomorphism,
$[x,a]^{\s}=[x,a]^{\tau}=[x^{\tau},a]=[x^{\s},a^{\s}],\forall a\in L_0,x\in L_1.$ Analogously,
$[x,a]^{\s}=[x^{\s},a^{\s}],\forall x\in L_{-1},a\in L_0.$ For $x,y\in L_1$ we have
$[(x+y),(x+y)^{\tau}]=[x,y^{\tau}]+[y,x^{\tau}]=[x,y^{\s}]+[y,x^{\s}]=0.$ Hence
$[x,y^{\tau}]^{\s}=[x,y^{\tau}]=[x^{\tau},y]=[x^{\s},y]$ and we proved that $\s$ is an automorphism
of $L.$ It is easy to see that a group generated by $\s,\r$ is the group $S_3$ and $L$ is a Lie algebra
with triality. $\Box$
\begin{Co}
Let $L$ be a Lie algebra with triality over a field $k$ as above. Then ${\cal F}(L)\simeq (L_1,\circ).$
\end{Co}
{\bf Proof.} Let $L=M\bo \L.$ Then for
$x=m\otimes (v-\xi w),$ $y=n\otimes (v-\xi w)$ by (\ref{f1}) we have:

$$
x\circ y=[x,y]^{\s}=(mn)_0\otimes (v-\xi w)^2_0+(mn)_2\otimes (v-\xi w)^2_2=(\xi-1)(mn)_2/3\otimes (v-\xi w),
$$
hence $x\circ y=(\xi-1)(mn)_2/3\otimes (v-\xi w)$ and $\phi:{\cal F}(L)=M\to L_1,\phi(m)=(\xi-1)m/3\otimes (v-\xi w)$
is an isomorphism. $\Box$
\begin{Ex}
Let $L$ be a split simple Lie algebra of type $D_4$ over a fild $k$ such that $1\neq\xi\in k,$
$\xi^3=1,$ and $\r_1$ be an automorphism
of $L$ of order 3 which admits triality (see Example 5). The group $AutL$ contains exactly
two conjugate classes of order 3: $\{\r_1\}$ and $\{\r_2\}$ (see (\cite{Kac}), Chapter 8).
We prove that $\r_2$ does not admit triality. We denote $\r_2=\r.$ Using (\cite{Kac}, Chapter 8)
we can write the bases of the spaces $L_0,$ $L_1$ and $L_{-1}.$

\[
\ba{llll}
L_0=&\{\xi\cd 1+2+\xi^2\cd 3, &\xi\cd 234+134+\xi^2\cd 124,&e_1+e_2+e_3,\\

    &f_1+f_2+f_3,&h_4,&h_1+h_2+h_3,\\

    &\xi^2\cd 12+13+\xi\cd 23,&\xi^2\cd 34+24+\xi\cd 14\},&
\ea
\]

\[
\ba{llll}

L_1=&\{4,           &14+24+34,&\xi\cd f_1+f_2+\xi^2\cd f_3,\\

    &1234,&124+134+234, &\xi^2\cd 1+2+\xi\cd 3, \\

&f_4,&\xi\cd h_1+h_2+\xi^2\cd h_3, &\xi\cd e_1+e_2+\xi^2\cd e_3,\\

   &&&\xi\cd 12+13+\xi^2\cd 23,\},
\ea
\]

\[
\ba{llll}

L_{-1}=&\{123,           &12+13+23,&\xi^2\cd f_1+f_2+\xi\cd f_3,\\

    &\emptyset,&1+2+3, &\xi^2\cd 14+24+\xi\cd 34, \\

&e_4,&\xi^2\cd h_1+h_2+\xi\cd h_3, &\xi^2\cd e_1+e_2+\xi\cd e_3,\\

   &&&\xi^2\cd 234+134+\xi\cd 124,\}.
\ea
\]

It is easy to see that $L_0$ is a split simple Lie algebra of type $A_2$ and $P$ is a Borel subalgebra with a basis: $\{\xi\cd 1+2+\xi^2\cd 3, e_1+e_2+e_3,
h_4,h_1+h_2+h_3,\xi^2\cd 12+13+\xi\cd 23\}.$ Moreover, $L_1$ and $L_{-1})$ are $L_0-$modules
with the highest weights $f_4 $ and $(123).$ But $[f_4,h_4]=-2f_4$ and $[123,h_4]=-(123).$ Hence $L_0-$modules
$L_1$ and $L_{-1}$ are not isomorphic. By Lemma 2 $\r$ does not admit triality.
\end{Ex}

Let $P=P_0\op P_2$ be a trivial Lie algebra with triality. We say that $P$ admits a $lifting$
if there exists a Lie algebra $L$ and an invariant subalgebra $Q$ in $T(L)$
such that $P/N(P)\simeq Q/N(Q).$ By definition, a Lie algebra $A\in {\cal T}$
admits a lifting if $\Phi(A)$ admits a lifting.
Let $A=\Psi(P)=A_0\op A_2$ and $N(P)=Ann_{A_0}A_2=0.$
Since $P$ is trival we have that $A$ is a Lie algebra and $A_0$ acts on the Lie algebra $A_2$
by derivations. Hence we have homomorphisms
$\phi: A_0\ra Der(A_2),$ $\psi: A_2\ra Inn(A_2)$ with $ker\phi=0,$ and
$ker\psi=C=\{x\in A_2\,|\, (xy)_2=0,\,\forall y\in A_2\}.$
Let $\ov{G}=\phi(A_0)+\psi(A_2)$ and $Q=\ov{G}/\psi(A_2).$ Then we have
the following short exact sequence
\bq
0\ra \psi(A_2)
{\buildrel\ov{i}\over\longrightarrow}\ov{G}{\buildrel\ov{j}\over\longrightarrow}Q\ra 0.
\label{xx}
\eq
Fix a section $s:Q\ra\ov{G}.$ Then for $x,y\in Q$
\bq
[s(x),s(y)]=s([x,y])+\ov{g}(x,y),
\label{xxx}
\eq
where $\ov{g}(x,y)\in \ov{i}(\psi(A_2)).$
Since $J(s(x),s(y),s(z))=0,$ for every $x,y,z\in Q$ we have by (\ref{xxx})
\bq
\ba{l}
\ov{g}([x,y],z)+\ov{g}(x,y)s(z)+\ov{g}([x,y],z)+\\[3mm]
\ov{g}(x,y)s(z)+\ov{g}([x,y],z)+\ov{g}(x,y)s(z)=0.
\ea
\label{y}
\eq
Let $t:\psi(A_2)\ra A_2$ be an arbitrary section. Then we can define
$g(x,y)=t(\ov{g}(x,y))\in A_2$ for $x,y\in Q.$
By (\ref{y}) we have that
\[
\ba{l}
f(x,y,z)=g([x,y],z)+g(x,y)s(z)+g([y,z],x)+\\[2mm]
g(y,z)s(x)+g([z,x],y)+g(z,x)s(y)\in C.
\ea
\]
Note that since $\ov{G}$ acts on $A_2,$ the product $A_2s(Q)$
is well defined and $C$ is a $Q$-module.
\begin{Prop}
Let $Z^3(Q,C)$ and $B^3(Q,C)$ be the groups of 3-cocycles and
3-coboundaries. Then with the above notation we have

(i)$\,f\in Z^3(Q,C),$

(ii)$\,f\in B^3(Q,C)$ if and only if $P$ admits lifting.
\end{Prop}
{\bf Proof.}
From the standard theory of extensions with non-Abelian
kernel \cite{Hoh},\cite{Mac} we have that
$f\in Z^3(Q,C)$ and $f\in B^3(Q,C)$ if and only if
the short exact sequence (\ref{xx})
may be lifted to the following short exact sequence
\bq
0\ra A_2
{\buildrel i\over\longrightarrow}G{\buildrel j\over\longrightarrow}Q\ra 0.
\label{z}
\eq
Suppose that $f\in B^3(Q,C)$ and we have the sequence (\ref{z}).
Then $(G,i(A_2))$ is a compatible pair of $G$ and for the corresponding
invariant subalgebra $B\subseteq T(G)$ we have $B/N(B)\simeq P/N(P).$
Conversely, if $P$ admits a lifting then there exists a Lie algebra $G$
and an invariant subalgebra $B\subseteq T(G)$ such that
$B/N(B)\simeq P/N(P).$ Hence holds for $G$ (\ref{z}).
$\Box$

%Recall that an algebra $L$ is almost algebraic if for $a\in L$
%there exists $b,c\in L$ such that $R_a=R_b+R_c,$ where $[R_b,R_c]=0,$ $R_b$ is
%a semisimple operator and $R_c$ is nilpotent.
 Denote by $Aut_{\li}L$ the group of $\li$-automorphisms of $L.$
In other words, $Aut_{\li}L=\{g\in Aut_k(L) |[g,S_3]=1\}.$

For any Malcev algebra $M$ and an ideal $I<M$ we define an ideal
${\cal G}_M(I)<{\cal G}(M)$ as follows
${\cal G}_M(I)=Inder(M,I)\oplus I,$ where $Inder(M,I)=\{D(x,y)=L_{[x,y]}+[L_x,L_y]\| x\in M,y\in I\}.$
It is clear that ${\cal G}_M(I)$ is a solvable ideal if and only if
$I$ is a solvable ideal.

\begin{Th}
Let $k$ be a field of characteristic 0 or $p>3$, $L$ be a finite dimensional algebra
in ${\cal L}(k),$ and
$M_2={\cal F}(L)\in {\cal M}(k).$

Then

(i) $M_2$ is semisimple (solvable) if and only if $L$ is semisimple (solvable).

(ii) $M_2$ is simple if and only if $L$ has no invariant ideals.

%(iii) $M_2$ is almost algebraic if and only if $L$ is.

(iii) $Aut_k(M_2)\simeq Aut_{\li}L.$

(iv) Let $G_2$ be the solvable radical of $M_2.$ Then there exists
a semisimple subalgebra $P_2$ such that $M_2=P_2\op G_2$ if and
only if there exists an invariant semisimple subalgebra $P$ in $L$
such that $L=P\op G,$ where $G$ is the solvable radical of $L.$

Moreover, for every two semisimple subalgebras $P_2$ and $Q_2$
such that $M_2=P_2\op G_2=Q_2\op G_2,$ there exists
$\phi\in Aut_k(M_2)$ such that $P_2^{\phi}=Q_2$
if and only if for every two semisimple invariant
subalgebras $P$ and $Q$ such that $L=P\op G=Q\op G$ there exists
$\psi\in Aut_{\li}L$
such that $P^{\psi}=Q.$

(v) Let $M_2$ be solvable and $N_2$ be the nilpotent radical of $M_2.$
Then there exists a torus $T_2$ such that $M_2=T_2\op N_2$
if and only if there exists an invariant torus $T$ in $L$ such
that $L=T\op N,$ where $N$ is the nilpotent radical of $L.$

Moreover, for every two tori $T_2$ and $Q_2$
such that $M_2=T_2\op N_2=Q_2\op N_2$ there exists
$\phi\in Aut_k(M_2)$ such that $T_2^{\phi}=Q_2$
if and only if for every two invariant
tori $T$ and $Q$ such that $L=T\op N=Q\op N$ there exists
$\psi\in Aut_{\li}L$
such that $T^{\psi}=Q.$

\end{Th}
{\bf Proof.}
(i) Let $M$ be a solvable Malcev algebra. Then $Inder(M)$ is a solvable
Lie algebra. Hence the Lie algebra with triality $L={\c G}(M)$ is solvable too.
If a Malcev algebra $M$ has a solvable ideal $I$ then
${\c G}_M(I)$ is a solvable ideal in $L.$ It follows from this
that the Malcev algebra $M$ is semisimple if $L$ is semisimple.
If $J$ is a solvable ideal of $L$ then $P=\sum_{g\in S_3}J^g$
is a solvable invariant ideal of $L$ and ${\c F}(P)$
is a solvable ideal of $M.$ If ${\c F}(P)=0$ then
$P\subseteq L_0$ and hence $P\subseteq Ann_{L_0}L_2=0.$

The other items of this Theorem are simple consequences of the definitions.
$\Box$

As a Corollary we can obtain the following Theorem.

\begin{Th}
Let  $M$ be a Malcev algebra over an algebraicly closed field $k$ of characteristic 0.

Then

(i) $M$ is semisimple if and only if $M$ is a direct sum
of Lie simple algebras and Malcev simple algebras of dimension 7.

(ii) Let $G$ be the solvable radical of $M.$ Then there exists
a semisimple subalgebra $P$ such that $M=P\op G.$

Moreover, for every two semisimple subalgebras $P$ and $Q$
such that $M=P\op G=Q\op G$ there exists
$\phi\in Aut_k(M)$ such that $P^{\phi}=Q.$

(iii)
Let $M$ be solvable. Then there exists a solvable almost algebraic Malcev algebra
 $R$ such that $M$ is a subalgebra of $R.$
If $R$ is an arbitrary almost algebraic solvable Malcev
algebra and $N$ is the nilpotent radical of $R$
then there exists a torus $T$ such that $R=T\op N.$

Moreover, for every two tori $T$ and $Q$
such that $R=T\op N=Q\op N$ there exists
$\phi\in Aut_k(R)$ such that $T^{\phi}=Q.$

\end{Th}
{\bf Proof.}
(i) Let $M$ be a semisimple Malcev algebra and $A={\c G}(M)$
be the corresponding Lie algebra with triality.
By Theorem 2, $A$ is semisimple and $A=A_1\op ... \op A_n.$
It is clear that for $g\in S_3,$ we have $A_i^g=A_{g(i)}.$
Hence every invariant minimal ideal $J$ of $A$
has the form $J=A_i\op A_j\op A_k$ or $J=A_i.$
In the first case we have the Lie algebra with triality from
Example 3. That is, that the corresponding
ideal in $M$ is a Lie direct summand.
We prove that in the second case $A_i$ is the simple Lie algebra
of type $D_4$ and $S_3$ is the group of diagram automorphisms.
Note that if $\l\in S_3$ with $\l^3=1,$  then $\l$ is not an inner automorphism.
Indeed, if $\l$ were an inner automorphism then there would exist a Cartan subalgebra
 $H$
in $A_i$ such that
$h^{\l}=h$ for $h\in H.$ But $A_i$ would be a Lie algebra with triality
hence $h^g=h,$ for all $g\in S_3$ and $h\in H.$
Let $A_i=\sum_{\a}\op A_{\a}$ be the Cartan decomposition of $A_i$
with $A_0=H.$
Then $A_{\a}=ke_{\a},\a\neq 0,$ and $e_{\a}^g=\b(g)e_{\a},\b(g)\in k,g\in S_3.$
 It follows from this that $S_3$ would be commutative, a contradiction.
However the group of outer automorphisms $Out(A_i)=Aut(A_i)/InAut(A_i)$
is the group of diagram automorphisms \cite{Jac} and the unique simple
Lie algebra with non-commutative group of outer automorphisms is an algebra
of type $D_4.$
In this case the corresponding Malcev algebra has dimension 7 (see Examples 5 and 6).

To prove item (ii) we need the following Lemma.
\begin{Le}\cite{Taft}
Let $S$ be a finite group of automorphisms of a finite dimensional Lie algebra $L$
and $|S|\neq 0\,(mod(char(k))),$ where $char(k)$ is the characteristic of the underlying
field $k.$

(i) If L has a Levi factor then $L$ has an $S$-invariant Levi factor.

(ii) If $L$ is a complete solvable Lie algebra, then there exists an $S$-invariant
torus $T$ such that $L=T\op N,$ where $N$ is the nilpotent radical of $L.$
\end{Le}
{\bf Proof.}
Let $P$ be a Levi factor of $L.$ Then $P\simeq L/G$ where $G$ is the radical of $L.$
Since $G$ is $S$-invariant, $S$ acts on $P.$
Fix an embedding $\phi: P\ra L.$
Define $\psi(p)=\sum_{g\in S}\phi(p^g)^{g^{-1}}.$ It is not difficult to
prove that $\psi$ is an embedding and $\psi(P)$ is an invariant
Levi factor.

Item (ii) of Lemma 3 can be proved in the same way.
$\Box$

Now item (ii) of Theorem 2 is a corollary of Lemma 3
and Taft's Theorem  (see \cite{Taft2}, Theorem 4).

Let $A$ be a solvable finite dimensional Lie algebra with triality and
$\tilde{A}$ be the completion of $A.$
It is clear that $S$ acts on $\tilde{A}$ and we prove that $\tilde{A}$ is an algebra with triality.
Denote by $N$ and $\tilde{N}$ the nil-radicals of $A$ and $\tilde{A}$
respectively and note that $N$ and $\tilde{N}$ are $S$-invariant. Hence,
by the construction of the completion we have an $S$-isomorphism
\bq
A/N\simeq \tilde{A}/\tilde{N}.
\label{p}
\eq
  It follows from Lemma 3 that
$\tilde{A}$ has an $S$-invariant torus $T$ such that $\tilde{A}=T\op \tilde{N}$
and it follows from (\ref{p}) that the $S$-module $T$ is of type $\L$.
From the construction of $\tilde{A}$ we have that for any $x\in \tilde{N}$ there
exists $t\in T$ such that $t-x\in A.$ Hence $\tilde{N}$ and $\tilde{A}$ are
 $S$-modules of type $\L.$
It is clear that the Malcev algebra ${\c G}(\tilde{A})$ is almost algebraic
and contains $M.$ Moreover, ${\c G}(\tilde{A})={\c G}(T)\op {\c G}(\tilde{N}).$
$\Box$

Note that this Theorem contains the main results of the following papers:
\cite{Carl1},\cite{Carl2},\cite{Gr1},\cite{Gr3},\cite{Kuz2},\cite{Sagl}.

As last application of the functors ${\cal F}$ and ${\cal G}$ to the theory
of Malcev algebras we prove that every Malcev algebra over a field
of characteristic $p>3$ may be embedded in a Malcev $p$-algebra.
 Note that the usual definition of Lie $p$-algebra (see, for example, \cite{Jac})
may be applied to a binary Lie algebra since for the definition we need only a fact that any two
generated subalgebra is a Lie algebra.
\begin{Prop}
Let $M$ be a Malcev algebra and $P$ be a Lie algebra with triality over a field $k$ of characteristic $p>3.$ Then
there exists a Malcev $p$-algebra $M_p$ and a Lie trial-$p$-algebra $P_p$
such that $M$ is an
ideal of $M_p$ and $P$ is an invariant ideal of $P_p.$ Moreover, $M^2\subseteq M_p$
and $P^2\subseteq P_p$.
\end{Prop}
{\bf Proof.}
The above statement for Malcev algebras is a consequence of the analogous statement
 for Lie algebras with triality and we have to prove only the Lie part of this Proposition.
Let $U(P)$ be the universal enveloping algebra of $P,$ and $Q$ be the $p$-closure of $P$ in $U(P).$
In other words $Q$ is the minimal $p$-subalgebra of Lie algebra $U(P)^{(-)}$ that contains $P$.
$Q$ has a basis $v_1, ... ,v_n$ where $v_i\in P$ or $v_i=x^{p^m}, \,x\in P.$ It is clear that every
automorphism $\phi$ of $P$ has extension $\phi_1\in Aut_kQ$. Hence $S_3$ acts on $Q,$ via:
$(x^{p^m})^{\phi}=(x^{\phi})^{p^m}.$
Thus we  have to prove that the $S_3$-module $Q/P$ has no one dimensional antisymmetric submodules.
If $kb$ were an one dimensional $S_3$-module, then $kb^{p}$ and $kb$ would be isomorphic $S_3$-modules.
If $kx\oplus ky$ were two dimensional irreducible $S_3$-module and $x^p,y^p$ are linearly independent
then $kx^p\oplus ky^p$ would be an irreducible $S_3$-module too.
Suppose that $y^p=\a x^p$ and $ x^{\s}=y,\, y^{\rho}=-x-y.$ Then $\a=\pm 1.$
On the
other hand $(y^p)^{\rho}=-x^p-y^p=(\a x^p)^{\rho}=\a y^p.$ Since the characteristic of the field $k$
is not 3, it follows that $x^p=y^p=0.$
  Now we can take every invariant ideal $I$ in $Q$ such that $I\cup P=0$
and put $P_p=Q/I.$
$\Box$
\begin{Co}
Let $M$ be a Malcev algebra over a field of characteristic $p>3.$
Then $M$ satisfies the following identity
$$
(xy.z)t^p+(yz.t^p)x+(zt^p.x)y-(xt^p.y)z=xz.yt^p,
$$
where $xt^p$ denotes $xR_t^p.$
\end{Co}

\section{Applications of the functors ${\c G}$ and ${\c F}$ to the theory
of Lie algebra with triality}

In this section we prove that in some sense every Lie algebra with triality
has an invariant ideal which is a trivial algebra with triality (see Definition 4)
and its factor algebra is an algebra with triality of type $D_4$ (see Example 4).

An algebra $A$ is perfect if $A^2=A.$
Every finite dimensional Lie algebra $A$ has a unique perfect ideal $A^{\omega}$
such that $A/A^{\omega}$ is a solvable algebra. It is clear that
$A^{\omega}=\cap_i A^{(i)}.$

\begin{Df}
Let $A=A_0\oplus A_2$ be a Lie algebra with triality and let
$T(A)_1=\{x\in A_2 \mid x^{\s}=x,\, \forall y\in A_2:
((xy)_2y)_2-((xy^{\s})_2y^{\s})_2-(x(yy^{\s})_2)_2=0\}.$
Then the space $T(A)=T(A)_0\oplus T(A)_2,$ where $T(A)_2=T(A)_1\oplus T(A)_1^{\rho},$
$T(A)_0=(T(A)_2A)_0,$ is called the T-centre of $A.$
\end{Df}
This definition seems to be absolutely artificial, but below we prove that it is analogous to
the notion of Lie centre in the theory of Malcev algebras.
Recall that the Lie centre of a Malcev algebra $M$ is the space
$L(M)=\{n\in M\mid J(n,x,y)=0, \forall x,y\in M\}.$
\begin{Prop}
Let $A$ be a Lie algebra with triality and $T(A)$ be the T-centre of $A$. Then ${\cal F}(T(A))$ is
the Lie centre of the Malcev algebra ${\cal F}(A).$
 Moreover, $T(A)$ is an invariant ideal of $A$ wich is trivial as Lie algebra with triality.
\end{Prop}
{\bf Proof.}
Let $A=A_0\oplus M\otimes (kv\oplus kw).$ Then for $x\in T_1(A)$ and $y\in M\otimes (kv\oplus kw)$
we have $x=n\otimes (v+w), \, y=p\otimes v+q\otimes w$ and
\[
\ba{l}
((xy)_2y)_2-((xy^{\s})_2y^{\s})_2-(x(yy^{\s})_2)_2=\\[4mm]
(np\cd p\otimes (v+w)-np\cd q\otimes w -nq\cd p\otimes v+nq\cd q\otimes (v+w))/9+\\[4mm]
(-nq\cd q\otimes (v+w)+nq\cd p\otimes w +np\cd q\otimes v-np\cd p\otimes (v+w))/9+\\[4mm]
(n(pq)\otimes w+n(qp)\otimes v)/9=\\[4mm]
\{((np)q+(pq)n+(qn)p)\otimes (v-w)\}/9=0.
\ea
\]
Hence $n\in L(M).$ Since $L(M)$ is an ideal of $M$ then
$T(A)={\cal G}_M(L(M))$ is an invariant ideal of $A.$
$\Box$

As a Corollary of Propositions 10 and 13 we have

\begin{Co}
Let $A$ be a normal Lie algebra with triality. Then $A$ is
trivial if and only if
$T(A)=A.$
\end{Co}

Recall that we denote by ${\c NG}$ the ${\bf Z}_2$-variety  defined by ${\bf Z}_2$-identities
(\ref{6}). In \cite{Gr5} we proved the following Theorem.
\begin{Th}
Let $\Gamma$ be the algebra of Example 2.
Then the Lie centre of a perfect finite  dimensional Malcev algebra $M$ over
 an algebraically closed field
of characteristic 0 is zero if and only if
$M\simeq B\bo\Gamma,$ where $B=B_0\op B_1\in {\c NG}$ and $B_0=B_1^2.$
\end{Th}

Let $B$ be a free algebra of the ${\bf Z}_2$-variety ${\c NG}$ with
odd generators $\{x_1,x_2, ...\}.$
Then $B=B_0\op B_1$ and $B_0=B_1^2.$
In \cite{Gr5} we proved the following Lemma.
\begin{Le}
The algebra $B_0$ is a commutative associative algebra with generators
$\{a_{ij}=-a_{ji}\mid 1\leq i<j\}$ and relations
\bq
a_{ij}a_{kp}+a_{jk}a_{ip}+a_{ki}a_{jp}=0.
\label{l}
\eq

Moreover, the proper words form a basis of $B$ where, by definition,
a word $w=a_{i_1j_1}a_{i_2j_2}\cd ... \cd a_{i_mj_m}\in B_0$ is proper
if

$
1)\, i_1\leq i_2\leq ... \leq i_m\leq j=min\{j_n\mid 1\leq n\leq m\};
$

\vspace{2mm}

2)\, If\, $i_a<i_b<i_c$ \,then\, $j_a<j_b<j_c$\, does\, not \,hold.

\vspace{2mm}

A word $w=a_{i_1j_1}a_{i_2j_2}\cd ... \cd a_{i_mj_m}x_k\in B_1$
is proper if

$
1)\, i_1\leq i_2\leq ... \leq i_m\leq j=min\{j_n\mid 1\leq n\leq m\};
$

\vspace{2mm}

2)\, If\, $k<i_m$ \,and \,$i_a<i_b<i_c$ \,then \,$j_a<j_b<j_c$\, does not hold;

\vspace{2mm}

3) \,If\, $i_m\leq k,$ \,then\, $k\leq j_m$ \,and

$
i_1\leq i_2\leq ... \leq i_m\leq j_p\leq ... \leq j_1,
k<j_{p+1}\leq ... \leq j_m,
$

where $i_{p+1}= ... =i_m=k$ and $i_p<k,p\leq m.$

\end{Le}
\begin{Co}
If $w,v\in B_0$ are proper words then the words
$wa_{1n}, vx_1$ and $vx_2$ are proper words too if
$w=a_{i_1j_1}a_{i_2j_2}\cd ... \cd a_{i_mj_m}$ and
$max\{j_p\mid 1\leq p\leq m\}\leq n.$
\end{Co}

\begin{Le} With the above notations we have

1. Every word on the generators $\{a_{ij}=-a_{ji}\mid 1\leq i<j\}$
is not a divisor of zero in $B_0.$

2. The generators $x_1$ and $x_2$ are linearly independent over $B_0.$
\end{Le}
{\bf Proof.}
1. It is enough to prove that any generator $a_{ij}$ is not a divisor of
zero in $B_0.$
If $a_{ij}$ were a divisor of zero in $B_0$ then $a_{ij}$ would be a divisor of zero
in some subalgebra $C_n$ of $B$ with finite generators $\{a_{ij}=-a_{ji}\mid 1\leq i<j\leq n\}.$
Let $v\in C_n$ and $va_{ij}=0.$
Suppose that $a_{ij}=a_{1,n}.$ We can order the set of proper words,
for instance, lexicographically. By Lemma 5, it follows that, for
arbitrary proper words
$w_1>w_2,$ we have that $w_1a_{1n}$ and $w_2a_{1n}$ are proper
words too. Moreover, $w_1a_{1n}>w_2a_{1n}.$
If $w$ is the maximal proper subword in $v$ then, from $va_{1n}=0,$ we have
that $wa_{1n}=0,$ a contradiction.
To finish the proof it is enough to note that any two arbitrary elements
$a_{ij}$ and $a_{pq}$ are conjugate in $Aut(B_0).$
Indeed, we have in the algebra $B:$ $a_{ij}=x_ix_j.$
As the elements $\{x_1,x_2, ...\}$ are a set of free generators of
$B$ then for any matrix $P=(\a_{ij})$ there exists
$\phi\in Aut(B)$
such that
$\phi(x_i)=\sum_j\a_{ij}x_j.$
It is clear that $B_0^{\phi}=B_0$ and some automorphism of this type
conjugates the elements $a_{ij}$ and $a_{pq}.$

2. Suppose that $vx_1+wx_2=0,$ for some $v,w\in B_0.$ Using Lemma 5
as above, we can obtain from this that $v_1x_1+w_1x_2=0,$ for some proper
words $v_1$ and $w_1.$ But this is impossible, since $v_1x_1$ and $w_1x_2$
are proper different words, by Corollary 4.
$\Box$

Denote by $R$ the $k$-algebra with generators $\{a_{ij},a_{ij}^{-1}\mid 1\leq i<j\}$
and relations (\ref{l}), and denote by $T\in {\c GN}$ the algebra which corresponds
to the simple 7-dimensional Malcev algebra. We have from \cite{Gr5}, that
%Note that $B\ot R\simeq T\ot R,$ where
$T=T_0\op T_1,T_0=ke,T_1=kx\op ky,$
$e^2=e,ex=x,ey=y,xy=e.$
\begin{Le}
With the above notations,
$$
R\ot_{B_0} B\simeq R\ot_kT.
$$
\end{Le}
{\bf Proof.} Define a morphism $\psi: R\ot_k T\ra R\ot_{B_0}B$ by the formula
$$\psi(\a e+\b x+\g y)=\a+\b x_1+\g a_{12}^{-1}x_2,$$
 where
$\a,\b,\g \in R.$
By Lemma 5, it follows that $ker\psi=0.$
By (\ref{l}), for $i>2,$ we have $x_i=x_1a_{2i}a_{12}^{-1}+x_2a_{1i}a_{12}^{-1}.$
Hence $Im\psi=R\ot_{B_0}B.$
$\Box$

Since $(R\ot_kT)\bo \Gamma\simeq R\ot_k(T\bo \Gamma)$ and
$(R\ot_{B_0}B)\bo\Gamma\simeq R\ot_{B_0}(B\bo\Gamma)$,
 by Lemma 6, we have
\bq
R\ot_{B_0}(B\bo\Gamma)\simeq R\ot_k(T\bo\Gamma).
\label{r}
\eq

Now we can prove the main result of this section.
\begin{Th}
Let $A$ be a finite dimensional perfect Lie algebra with triality over an algebraically
closed field $k$ of characteristic 0.

Then, with the above notations, we have

1. $L(A/L(A))=0.$

2. There exists an $S$-invariant $B_0$-subalgebra $F$ in the Lie algebra with triality
$R\ot_k D_4$ such that $A$ is a homomorphic image of $F$ if $L(A)=0.$ Moreover,
$F$ does not depend on $A.$
\end{Th}
{\bf Proof.}
Let $M_2={\c F}(A)$ be the corresponding Malcev algebra. Since $A$ is perfect,
 $M_2$ is perfect too. By Lemma 5 \cite{Gr5}, it follows that
$L(M_2/L(M_2))=0;$ hence $L(A/L(A))=0.$

Suppose that $L(A)=L(M_2)=0.$ Then by Proposition 1 \cite{Gr5} one gets that
$M_2$ is a homomorphic image of $B\bo\Gamma .$ By (\ref{r})
we have that $B\bo\Gamma$ is a $B_0$-subalgebra of $R\ot_k(T\bo\Gamma).$
But this Malcev algebra corresponds to a Lie algebra of the type $D_4$
over $R$ (see Example 5.)
$\Box$

Let $M$ be the free Malcev algebra with a countable set of free generators.
Then $F={\cal G}(M)$ is a free Lie algebra with triality in the following sense.
For every normal Lie algebra with triality $L=L_0\oplus L_2$ there exists
an epimorphism ${\cal G}(M)\ra L.$ We can construct $F$ as follows.
Let $P$ be the free Lie algebra with free generators
$\{x_1,y_1,x_2,y_2, \,...\}.$
The group $S_3$ acts on $P$ by the rule
$x_i^{\s}=y_i,\,x_i^{\r}=y_i,\,y_i^{\r}=-x_i-y_i,i=1,2,...$ .
Then $P=P_0\oplus P_1\oplus P_2,$
where $P_0=\{x\in P\, |\, x^{\l}=x,\,\forall \l\in S_3\},$
$P_1=\{x\in P\, |\, x^{\s}=-x,\,x^{\r}=x\},$
and $P_2$ is the sum of 2-dimensional irreducible $S_3$-modules.
We denote by $I$ the ideal of $P$ generated by $P_1.$
Then $F=P/I.$

In the series of papers \cite{Fi1}-\cite{Fi8}, V. Filippov developed the deep theory
of free Malcev algebras. As a rule it is not difficult to reduce the results
of Filippov about Malcev algebras to Theorems about Lie algebra with triality.
For example, Filippov proved \cite{Fi6} that $Z(M)=\{x\in M\,|\, xM=0 \} \neq 0.$
Hence we have that $Z(F)\neq 0.$

From the main result of \cite{Fi3} we are obtain:

\begin{Th}
Let $A$ be a simple central Lie algebra with triality over an algebraicly
closed field of characteristic $\neq 2,3.$
Then $A$ is a finite dimensional simple Lie algebra of type $D_4.$
\end{Th}

Let $A$ be a Lie algebra with triality and $V$ be an $A$-module. Then by definition
$V$ is an $A$-module with triality if the extension $A\oplus V$ is a Lie algebra with triality,
where $V$ is an Abelian invariant ideal.
It is clear that if $V$ is an $A$-module with triality then
${\cal F}(V)$ is a Malcev module over ${\cal F}(A)$
and conversely. There exists an example of finite dimensional nilpotent Malcev algebra $Q$ \cite{Fi8}
such that $Q$ has no exact finite dimensional Malcev modules.
Hence we have an example of finite dimensional nilpotent Lie algebra with triality ${\cal G}(Q)$
without exact finite dimensional modules with triality.

I wish to express my gratitude to Professor A. Glass for the help
in preparation of this paper and to referee for many usefull remarks that improve the article,

\end{document}